\documentclass[a4paper, 12pt]{article}
\usepackage{amsfonts}
\usepackage{fancyhdr}
\usepackage{graphicx}
\usepackage{epsfig}
\usepackage{amssymb}
\usepackage{mathrsfs,amsfonts,amsmath}
\usepackage{color}
\usepackage{amsmath}
\usepackage{caption}
\usepackage{graphicx, subfig}
\allowdisplaybreaks[4]
\makeatletter

\@addtoreset{equation}{section} \makeatother
\newtheorem{Assumption}{Assumption}[section]
\newtheorem{Theorem}[Assumption]{Theorem}

\newtheorem{Remark}[Assumption]{Remark}
\newtheorem{Lemma}[Assumption]{Lemma}

\begin{document}
	\title{\textbf{Mean square exponential stability of numerical methods for stochastic differential delay equations
			\footnote{Supported by Natural Science Foundation of Beijing Municipality (1192013).}}}
	\author{
		Guangqiang Lan\footnote{Corresponding author: Guangqiang Lan:
			langq@buct.edu.cn.} \quad and\quad Qi Liu\footnote{Qi Liu: 2021201048@buct.edu.cn.}
		\\ \small College of Mathematics and Physics, Beijing University of Chemical Technology, Beijing 100029, China}
	\date{}	
	\maketitle	
	\begin{abstract}
	Mean square exponential stability of $\theta$-EM and modified truncated Euler-Maruyama (MTEM) methods for stochastic differential delay equations (SDDEs) are investigated in this paper. We present new criterion of mean square exponential stability of the $\theta$-EM and MTEM methods for SDDEs, which are different from most existing results under Khasminskii-type conditions. Two examples are provided to support our conclusions.
		
	\end{abstract}
	
	
	
	\section{Introduction}
	
	\noindent Consider the nonlinear SDDE
	\begin{equation}\label{sdde}
		d x(t)=f(x(t), x(t-\tau(t))) dt+g(x(t), x(t-\tau(t)))dB(t), t\ge0
	\end{equation}
	with initial value $\xi$, $B(t):=\left(B_{1}(t), \ldots,B_{m}(t)\right)^{T}$ be an $m$-dimensional standard Brownian motion, $f: \mathbb{R}^{d} \times \mathbb{R}^{d} \rightarrow \mathbb{R}^{d}, \quad g: \mathbb{R}^{d} \times \mathbb{R}^{d}  \rightarrow \mathbb{R}^{d \times m}$ and $\tau(t): \mathbb{R}_{+} \rightarrow \mathbb{R}_+$ are Borel-measurable functions and $0<\tau(t) \leq \tau, t \in \mathbb{R}_{+}$, for some $\tau>0$.

Asymptotic stability of numerical approximations for the above SDDE (\ref{sdde}) or more general model has been widely investigated in recent years. In general, to obtain the exponential stability results, the following Khasminskii-type conditions are usually parts of the sufficient conditions.
\begin{equation}\label{K}
		2\langle x,f(x,y)\rangle+||g(x,y)||^2\le -C_1|x|^2+C_2|y|^2
	\end{equation}
where $C_1>C_2>0.$

For example, in \cite{14}, Chapter 5, Razumikhin type Theorems are presented for SDDEs, where conditions in Corollary 6.6 implies \eqref{K} if $k=1$. In \cite{1}, Theorem 4.1 implied that if there is no neutral term, then the MTEM method is mean square exponentially stable under condition (\ref{K}). Under the same condition, \cite{RAM} considered exponential mean-square stability of two classes of theta Milstein methods for stochastic differential equations with constant delay. \cite{ZSL} obtained exponential stability of stochastic theta method for nonlinear stochastic differential equations with piecewise continuous arguments (in this case, $\tau(t)=t-[t]$). \cite{LMD} obtained mean square exponential stability of split-step method under stronger conditions (note that they need the linearity of both $f$ and $g$), while \cite{LZ} and \cite{OM} investigated mean square stability of two class of theta method of neutral stochastic differential delay equations under similar conditions.

Recently, by considering each component separately, \cite{NH} presented a different type of sufficient conditions under which the trivial solution of given SDDE is mean square exponentially stable.

Motivated by that paper, now suppose that $f$ and $g$ jointly satisfy the following
\begin{equation}\label{weak}
2x_{i} f_{i}(x, y) +\sum_{l=1}^{m}\left(g_{il}(x,y)\right)^{2}\leq \sum_{j=1}^{d}a_{i j}x_{j}^{2}+\sum_{j=1}^{d} b_{i j}y_{j}^{2}, \quad i\in\underline{d}:={1,2,...,d}.
\end{equation}

By the same method in \cite{NH}, it is not difficult to obtain mean square exponential stability of $x(t)$ under \eqref{weak} if $a_{ij}$ and $b_{ij}$ satisfy suitable conditions. Note that $f$ and $g$ in \eqref{weak} do not necessarily satisfy linear growth condition.

Note also that by (\ref{weak}) we can only have
$$2\langle x,f(x,y)\rangle+||g(x,y)||^2\le \sum_{j=1}^{d}\sum_{i=1}^{d}a_{i j}x_{j}^{2}+\sum_{j=1}^{d} \sum_{i=1}^{d}b_{i j}y_{j}^{2}.$$

If $\sum_{i=1}^{d}a_{i j}\ge0$ for some $j\in\underline{d}$, or $0<\min_{j\in\underline{d}}(-\sum_{i=1}^{d}a_{i j})<\max_{j\in\underline{d}}\sum_{i=1}^{d}b_{i j}$, then (\ref{K}) can never hold. Therefore, (\ref{weak}) is a different type condition from (\ref{K}).

To the best of our knowledge, there isn't any result of stability of numerical methods under condition (\ref{weak}). So the key aim of this paper is to obtain mean square exponential stability of $\theta$-EM method and MTEM method of the given stochastic differential delay equations under (\ref{weak}) with suitable $a_{ij}$ and $b_{ij}$.

The rest of the paper is organized as follows. In Section 2, $\theta$-EM method, MTEM method and some necessary assumptions will be presented. In Section 3, we will prove that the given $\theta$-EM method is mean square exponentially stable under (\ref{weak}). Section 4 will prove that the MTEM method will replicate the mean square exponential stability of the exact solution for the stochastic differential delay equation under the same conditions. Finally, in Section 5, two examples will be presented to interpret our conclusion.

\section{Model description and preliminaries}
	
Let $\mathbb{N}$ be the set of all natural numbers. Denote by $\mathbb{R}$ the set of all real numbers. For given integers $l, q \geq 1, \mathbb{R}^{l}$ is the $l$-dimensional vector space over $\mathbb{R}$ and $\mathbb{R}^{l \times q}$ stands for the set of all $l \times q$ matrices with entries in $\mathbb{R}$.  For given $m \in \mathbb{N}$, let $\underline{m}:=\{1,2, \ldots, m\}$.
	Let $\left(\Omega, \mathcal{F},\{\mathcal{F}\}_{t \geq 0}, P\right)$ be a complete probability space with a filtration $\left\{\mathcal{F}_{t}\right\}_{t \geq 0}$ satisfying the usual conditions (i.e., right continuous and $\mathcal{F}_{0}$ containing all $P$ null sets).  Additionally, $\|\cdot\|$ denotes the Hilbert-Schmidt norm $\|g\|^2:=\sum_{i=1}^d\sum_{j=1}^mg^2_{ij}$ for any $d\times m$ matrix $g=(g_{ij})\in\mathbb{R}^{d\times m}$. Let $C\left([-\tau, 0], \mathbb{R}^{d}\right)$ be the Banach space of all continuous functions on $[-\tau, 0]$, endowed with the norm $\|\varphi\|=$ $\max _{s \in[-\tau, 0]}\|\varphi(s)\|$.
	Denote by $C_{\mathcal{F}_0}^{b}\left([-\tau, 0], \mathbb{R}^{d}\right)$, the family of $\mathcal{F}_0$-measurable bounded $C\left([-\tau, 0], \mathbb{R}^{d}\right)$-valued random variables.
	
	An initial condition for (\ref{sdde}) is defined by
	\begin{equation}\label{2}
		x\left(s\right)=\xi(s), \quad s \in[-\tau, 0]
	\end{equation}
	where $\xi \in C_{\mathcal{F}_0}^{b}\left([-\tau, 0], \mathbb{R}^{d}\right)$.
	
	An $\mathbb{R}^{d}$-valued stochastic process $\{x(t)\}_{t \geq -\tau}$ is called a solution of (\ref{sdde}) with initial value (\ref{2}), if $\{x(t)\}_{t \geq 0}$ is continuous, $\mathcal{F}_{t}$-adapted such that
	\begin{equation}\label{3}
		\begin{aligned}
			x(t)= \xi(0)+\int_0^{t} f(x(s), x(s-\tau(s))) d s
			+\int_0^{t} g(x(s), x(s-\tau(s))) d B(s)
		\end{aligned}
	\end{equation}
	holds for each $t \geq 0$, with probability 1 and the initial condition (\ref{2}) is fulfilled.

	\begin{Assumption}\label{assumption 2.1}
		Assume that both the coefficients $f$ and $g$ in (\ref{sdde}) are locally Lipschitz continuous, that is, for each $R>0$ there is $L_{R}>0$ $($depending on $R$$)$ such that
		\begin{equation}\label{2.4}
			|f( x, y)-f(\bar{x}, \bar{y})| \vee|g(x, y)-g(\bar{x}, \bar{y})| \leq L_{R}(|x-\bar{x}|+|y-\bar{y}|)
		\end{equation}
		for all $|x| \vee|y| \vee|\bar{x}| \vee|\bar{y}| \leq R>0$.
	\end{Assumption}
	
It is obvious that $L_{R} $ is an increasing function with respect to $R $.

Let $\Delta$ be a stepsize such that $\tau=\bar{m}\Delta$ for some positive integer $\bar{m}$. Then for any $\theta\in[0,1],$ we can define $\theta$ Euler-Maruyama method ($\theta$-EM for short) $X_{k}$ as the following:
	\begin{equation}\label{theta}
		\begin{aligned}
			&X_{k}=\xi(k \Delta), \quad k=-\bar{m},-\bar{m}+1, \ldots, 0 .\\
&X_{k+1}=X_{k}+\left((1-\theta)f\left(X_{k}, X_{k-[\frac{\tau(k\Delta)}{\Delta}]}\right)+\theta f\left(X_{k+1}, X_{k+1-[\frac{\tau((k+1)\Delta)}{\Delta}]}\right)\right) \Delta\\&
			\qquad+g\left( X_{k}, X_{k-[\frac{\tau(k\Delta)}{\Delta}]}\right) \Delta B_{k}, \quad k=0,1,2, \ldots,\\
		\end{aligned}
	\end{equation}
	Here $\Delta B_{k}=B((k+1) \Delta)-B(k \Delta)$ is the increment of the $m$-dimensional standard Brownian motion.

Note that if $\theta=0$, it becomes to classical EM method, if $\theta=1$, it is called Backward Euler method. Moreover, since it is an implicit method for $\theta\in(0,1]$, then it is necessary to make sure that $\theta$-EM method is well defined. So we need the following assumption

	\begin{Assumption}\label{danbian}
		Assume that $f$ is one-sided Lipschitz continuous, that is, there is $L>0$ such that
		\begin{equation}\label{db}
			\left\langle x-\bar{x}, f(x,y)-f(\bar{x},y)\right\rangle  \leq L|x-\bar{x}|^2
		\end{equation}
		for all $x,y,\bar{x}\in\mathbb{R}^d$.
	\end{Assumption}

According to \cite{WMS}, if $L\theta\Delta<1$, then the $\theta$-EM scheme \eqref{theta} is well defined, see also \cite{HW,MY} for details.

Now let us give the modified truncated Euler-Maruyama method for (\ref{sdde}).
	
For  $\Delta^{*}>0 $ ,  let  $ h(\Delta) $ be a strictly positive decreasing function $h:\left(0, \Delta^{*}\right] \rightarrow(0, \infty)$ such that
	$$
	\lim _{\Delta \rightarrow 0} h(\Delta)=\infty \text { and } \lim _{\Delta \rightarrow 0} L_{h(\Delta)}^{2} \Delta=0 .
	$$

According to \cite{LX}, Remark 2.1, such $h$ always exists.

	For any $\Delta \in\left(0, \Delta^{*}\right)$, we define the modified truncated function of $f$ as the following:
	$$
	f_{\Delta}( x, y)=\left\{\begin{array}{lc}
		f( x, y),\quad   \quad  \quad   \quad  \quad  \quad \quad |x| \vee|y| \leq h(\Delta), \\
		\frac{|x| \vee|y|}{h(\Delta)} f\left(\frac{h(\Delta)}{|x| \vee|y|}(x, y)\right),\quad |x| \vee|y|>h(\Delta) .
	\end{array}\right.
	$$
	
	$g_{\Delta}$ is defined in the same way as $f_{\Delta}$. Here $f(a(x, y)) \equiv f(a x, a y)$ for any $t \geq 0, a\in(0,1), x, y \in \mathbb{R}^{d}$.
	It is obvious that the functions $f_{\Delta}$ and $g_{\Delta}$ defined above are unbounded while the truncated functions defined in $\cite{8,9}$ are bound for any fixed $\Delta$.
	
Then by using $f_{\Delta}$ and $g_{\Delta}$, we can define the modified truncated Euler-Maruyama method (MTEM for short) $X_{k}$ as the following:
	\begin{equation}\label{em}
		\begin{aligned}
			&X_{k+1}=X_{k}+f_{\Delta}\left( X_{k}, X_{k-[\frac{\tau(k\Delta)}{\Delta}]}\right) \Delta
			+g_{\Delta}\left( X_{k}, X_{k-[\frac{\tau(k\Delta)}{\Delta}]}\right) \Delta B_{k}, \quad k=0,1,2, \ldots, \\
			&X_{k}=\xi(k \Delta), \quad k=-\bar{m},-\bar{m}+1, \ldots, 0 .\\
		\end{aligned}
	\end{equation}

\begin{Assumption}\label{assumption 2.2}
		
Let $f(x, y):=\left(f_{1}(x, y), \ldots, f_{d}(x, y)\right)^{T} \in \mathbb{R}^{d}$ and $g(x,y):=\left(g_{i l}(x, y)\right) \in \mathbb{R}^{d \times m}$. Suppose there exist constants $a_{ii}\in \mathbb{R};a_{ij}\ge0,i\neq j ;b_{ij}\geqslant0,$ $i,j  \in \underline{d}$, such that
 \begin{equation}\label{w}
2x_{i} f_{i}(x, y) +\sum_{l=1}^{m}\left(g_{i l}(x,y)\right)^{2}\leq \sum_{j=1}^{d}a_{i j}x_{j}^{2}+\sum_{j=1}^{d} b_{i j}y_{j}^{2}
\end{equation}
 holds for any $i\in \underline{d}$.
	\end{Assumption}

It follows easily from Assumption \ref{assumption 2.2} that
	\begin{equation}\label{3.1}
		2\left\langle x,f(x,y)\right\rangle +|g(x,y)|^{2}\leq K(1+|x|^{2}+|y|^{2})	
	\end{equation}
	where $K=\max\limits_{j}\left\lbrace \sum\limits_{i=1}^d|a_{ij}|,\sum\limits_{i=1}^db_{ij}\right\rbrace\vee 1$.

Similar to Theorem 3.5 in\cite{GMY} or Theorem 2.1 in \cite{LW}, it is easy to see that Assumptions \ref{assumption 2.1} and \ref{assumption 2.2} guarantee strong convergence of MTEM method. If Assumption \ref{danbian} also holds, then, by \cite{ZYY}, the $\theta$-EM method is also strongly convergent.
	
The following Theorem provides explicit criteria for the mean square exponential stability of (\ref{sdde}).

\begin{Theorem}\label{Theorem1}
Suppose Assumptions \ref{assumption 2.1} and \ref{assumption 2.2} hold. If there exist constants $p_{j}>0,j\in \underline{d}$ such that
		\begin{equation}\label{changshu}
			\sum_{j=1}^{d}(a_{i j}+b_{i j})p_{j}< 0,
		\end{equation}
then the initial value problem (\ref{sdde}), (\ref{2}) always has a unique solution $x\left(t, \xi\right)$ and the trivial solution of (\ref{sdde}) is exponentially stable in mean square sense.
\end{Theorem}

Note that \eqref{changshu} implies $\max_{i\in \underline{d}}a_{ii}<0$ since $a_{ij}\ge0,i\neq j ;b_{ij}\geqslant0,$ $i,j  \in \underline{d}$.

\textbf{Proof of Theorem \ref{Theorem1}.} Firstly, local Lipschitz continuity of $f$ and $g$ and (\ref{3.1}) implies that (\ref{sdde}) always has a unique local maximum solution. Moreover, since Assumption \ref{assumption 2.2} implies \eqref{3.1}, then the solution must be non explosive, that is, \eqref{sdde} has a unique global solution $x\left(t, \xi\right)$.
	
 Let $p:=$ $\left(p_{1}, p_{2}, \ldots, p_{d}\right)^{T} \in \mathbb{R}^{d}$. By continuity, \eqref{changshu} implies that there exists $\beta>0$ small enough such that
	\begin{equation}\label{9}
		\sum_{j=1}^{d}\left(a_{i j}+e^{\beta \tau} b_{i j}\right) p_{j} \leq- \beta p_{i}
	\end{equation}
	for any $i \in \underline{d}$.
	
	Then by the same method of the proof of Theorem II.2, \cite{NH}, it follows that the trivial solution of \eqref{sdde} is exponentially stable in mean square sense. $\square$

\section{Mean exponential stability of $\theta$-EM method}

In this Section, we consider mean exponential stability of $\theta$-EM method for all $\theta\in[0,1]$.

\begin{Theorem}\label{Theorem21}
		Suppose Assumption \ref{assumption 2.1}, Assumption \ref{danbian}  and Assumption \ref{assumption 2.2} hold. Moreover, there exist constants $0<\varepsilon<\frac{\min_i|a_{ii}|}{d\max_i|a_{ii}|}$ and $p_{j}>0,j\in \underline{d}$ such that
		\begin{equation}\label{changshu3}
			\varepsilon a_{ii}p_i+\sum_{j\neq i}a_{i j}p_{j}+\sum_{j=1}^{d}b_{i j}p_{j}< 0.
		\end{equation} Then for any fixed $\theta\in(\frac{1}{2},1]$, the $\theta$-EM method \eqref{theta} is mean quare and almost surely exponentially stable.
	\end{Theorem}	
	
\begin{Remark}\label{R1}
Notice that \eqref{changshu3} implies
\begin{equation}\label{changshu2}
			\sum_{j=1}^da_{i j}p_{j}+\sum_{j=1}^{d}b_{i j}p_{j}< 0.
		\end{equation}
Then similar to the proof of Theorem II.2 \cite{NH}, it follows that the Assumption \ref{assumption 2.2} implies that the solution $x\left(t, \xi\right)$ to (\ref{sdde}) is exponentially stable in mean square sense if \eqref{changshu} holds. Theorem \ref{Theorem21} assures that $\theta$-EM method \eqref{theta} replicates mean square exponential stability of the exact solution under given conditions.
\end{Remark}

	\textbf{Proof of Theorem \ref{Theorem21}.}
	By \eqref{theta}, for any $k\geqslant 0$ and $i=1,2,...,d$, we have
	\begin{equation}\label{3.70}
		\begin{aligned}
			F_{k+1}^i&=F_k^i+f_{i}\left( X_{k}, X_{k-[\frac{\tau(k\Delta)}{\Delta}]}\right) \Delta+\sum_{l=1}^{m}g_{il}\left( X_{k}, X_{k-[\frac{\tau(k\Delta)}{\Delta}]}\right) \Delta B_{k}
		\end{aligned}
	\end{equation}
where $F_k^i=X_{k}^{i}-\theta f_{i}\left( X_{k}, X_{k-[\frac{\tau(k\Delta)}{\Delta}]}\right) \Delta$.

So,
	\begin{equation}\label{3.80}
		\begin{aligned}
			|F_{k+1}^{i}|^{2}&=|F_{k}^{i}|^{2}+\left(2X_{k}^{i} f_{i}\left( X_{k}, X_{k-[\frac{\tau(k\Delta)}{\Delta}]}\right)
			+\sum_{l=1}^{m}g_{il}^{2}\left( X_{k}, X_{k-[\frac{\tau(k\Delta)}{\Delta}]}\right) \right)\Delta\\&\quad +(1-2\theta)f^{2}_{i}\left(X_{k}, X_{k-[\frac{\tau(k\Delta)}{\Delta}]}\right) \Delta^{2} +M^i_{k}\\
		\end{aligned}
	\end{equation}
	where
	$$\begin{aligned}
		M^i_{k}&=2\left(F_k^i+ f_{i}\left( X_{k},X_{k-[\frac{\tau(k\Delta)}{\Delta}]}\right)\Delta\right)\sum_{l=1}^{m}g_{il}\left( X_{k}, X_{k-[\frac{\tau(k\Delta)}{\Delta}]}\right)\Delta B^l_{k}\rangle\\
		&\quad+\left(\left|\sum_{l=1}^{m}g_{il}\left( X_{k}, X_{k-[\frac{\tau(k\Delta)}{\Delta}]}\right)\Delta B^l_{k}\right|^{2}-\sum_{l=1}^{m}g^2_{il}\left( X_{k}, X_{k-[\frac{\tau(k\Delta)}{\Delta}]}\right)\Delta\right).
	\end{aligned}$$

	It is easy to see that $M^i_{k}$ is a $\mathscr{F}_{k\Delta}$ martingale and $ EM^i_{k}=0$. Moreover, we claim that
\begin{equation}\label{claim}\aligned &\quad2X_{k}^{i} f_{i}( X_{k},X_{k-[\frac{\tau(k\Delta)}{\Delta}]})+\sum_{l=1}^mg^2_{il}( X_{k},X_{k-[\frac{\tau(k\Delta)}{\Delta}]})+(1-2\theta)f^{2}_{i}\left(X_{k}, X_{k-[\frac{\tau(k\Delta)}{\Delta}]}\right) \Delta\\&
\le -C|F_k^i|^2+\sum_{j\neq i}^{d}a_{i j}|X_{k}^j|^{2}+\sum_{j=1}^{d} b_{i j}|X_{k-[\frac{\tau(k\Delta)}{\Delta}]}^j|^{2}\endaligned\end{equation}
if $C<\min_{i}|a_{ii}|.$

In fact, denote $f_k^i:=f_{i}\left(X_{k}, X_{k-[\frac{\tau(k\Delta)}{\Delta}]}\right), g_k^{il}:=g_{il}\left(X_{k}, X_{k-[\frac{\tau(k\Delta)}{\Delta}]}\right)$ it follows that
$$\aligned (2\theta-1)|f_k^i|^2\Delta-C|F_k^i|^2&=((2\theta-1)\Delta-C\theta^2\Delta^2)|f_k^i|^2+2C\theta\Delta X_{k}^{i} f^i_{k}-C|X_{k}^{i}|^2\\&
=a(f_k^i+bX_k^i)^2-(ab^2+C)|X_k^i|^2\endaligned$$
where $a=(2\theta-1)\Delta-C\theta^2\Delta^2,$ $b=\frac{C\theta\Delta}{a}.$

Now we can choose $\Delta$ small enough such that $a>0$ since $\theta>\frac{1}{2}$. Moreover, $ab^2+C\le -a_{ii}, \forall i=1,2,...,d$, and $a_{ii}<0$ by \eqref{changshu}. Then by Assumption \ref{assumption 2.2}, we have

$$\aligned (2\theta-1)|f_k^i|^2\Delta-C|F_k^i|^2&\ge-\frac{ab^2+C}{a_{ii}}a_{ii}|X_k^i|^2\\&\ge a_{ii}|X_k^i|^2\\&
\ge 2X_{k}^{i} f_k^{i}+\sum_{l=1}^m(g_k^{il})^2-(\sum_{j\neq i}^{d}a_{i j}|X_{k}^j|^{2}+\sum_{j=1}^{d} b_{i j}|X_{k-[\frac{\tau(k\Delta)}{\Delta}]}^j|^{2}) \endaligned$$
which implies \eqref{claim}.

Then by \eqref{3.80} and \eqref{claim}, it follows that
$$|F_{k+1}^{i}|^{2}\le (1-C\Delta)|F_k^i|^2+(\sum_{j\neq i}^{d}a_{i j}|X_{k}^j|^{2}+\sum_{j=1}^{d} b_{i j}|X_{k-[\frac{\tau(k\Delta)}{\Delta}]}^j|^{2})\Delta+M^i_k.$$

Therefore for any $\alpha>0$,
$$\aligned |F_{k+1}^{i}|^{2}e^{\alpha(k+1)\Delta}-|F_k^i|^2e^{\alpha k\Delta}&\le (1-C\Delta-e^{-\alpha\Delta})|F_k^i|^2e^{\alpha(k+1)\Delta}+(\sum_{j\neq i}^{d}a_{ij}e^{\alpha(k+1)\Delta}|X_{k}^j|^{2}\\&\quad+\sum_{j=1}^{d} b_{i j}e^{\alpha(k+1)\Delta}|X_{k-[\frac{\tau(k\Delta)}{\Delta}]}^j|^{2})\Delta+M^i_ke^{\alpha(k+1)\Delta}.\endaligned$$

Then
$$\begin{aligned}
		 e^{\alpha k\Delta}E|F_{k}^{i}|^{2}&\le E|F_{0}^{i}|^{2}+R(\Delta)\sum_{l=0}^{k-1}e^{\alpha(l+1)\Delta}E|F_{l}^{i}|^{2}\\&
\quad+(\sum_{j\neq i}^{d}a_{ij}\sum_{l=1}^{k-1}e^{\alpha(l+1)\Delta}E|X_{l}^{j}|^{2}
+\sum_{j=1}^{d}b_{ij}\sum_{l=1}^{k-1}e^{\alpha(l+1)\Delta}E|X_{l-[\frac{\tau(l\Delta)}{\Delta}]}^{j}|^{2})\Delta
	\end{aligned}$$
where $R(\Delta)=1-C\Delta-e^{-\alpha\Delta}.$

Since $R(0)=0$ and
$$R'(\Delta)=-C+\alpha e^{\alpha\Delta}\le0$$
for all $\Delta>0$ if $\alpha\le C$, then $R(\Delta)\le0.$

It is obvious that if we choose $C'>E||\xi||^{2}\max_i \frac{1}{p_i}\vee \max_i\frac{E|F_{0}^{i}|^{2}}{\varepsilon' p_i}$, where $\varepsilon'$ will be given in the following, then we have
	\begin{equation}\label{chushi1}
		E|X_{k}^{i}|^{2}\leq C'p_{i}e^{-\beta k\Delta} ,\quad k=-\bar{m},-\bar{m}+1, \ldots, 0, i=1,2,...,d.
	\end{equation}

Now we will prove that (\ref{chushi1}) also holds for any $k>0$ by using second principle of mathematical induction.

Assume that
\begin{equation}\label{jiashe1}E|X_{l}^{i}|^{2}\le C'p_ie^{-\beta l\Delta}, \forall -\bar{m}\le l\le k-1, i=1,2,...,d.\end{equation}

Then
$$\begin{aligned}
		 e^{\alpha k\Delta}E|F_{k}^{i}|^{2}&\le E|F_{0}^{i}|^{2}+\left(\sum_{j\neq i}^{d}a_{ij}\sum_{l=1}^{k-1}e^{\alpha(l+1)\Delta}E|X_{l}^{j}|^{2}
+\sum_{j=1}^{d}b_{ij}\sum_{l=1}^{k-1}e^{\alpha(l+1)\Delta}E|X_{l-[\frac{\tau(l\Delta)}{\Delta}]}^{j}|^{2}\right)\Delta\\&
\le|F_{0}^{i}|^{2}+\Delta\sum_{j\neq i}^{d}a_{ij}\sum_{l=1}^{k-1}e^{\alpha(l+1)\Delta}C'p_je^{-\beta l\Delta}\\&
\quad+\Delta\sum_{j=1}^{d}b_{ij}\sum_{l=1}^{k-1}e^{\alpha(l+1)\Delta}C'p_je^{-\beta(l-[\frac{\tau(l\Delta)}{\Delta}])\Delta}\\&
\le |F_{0}^{i}|^{2}+C'\Delta\left(\sum_{j\neq i}^{d}a_{ij}p_j+\sum_{j=1}^{d}b_{ij}p_je^{\beta\tau}\right)e^{\alpha\Delta}\frac{e^{(\alpha-\beta)k\Delta}-1}{e^{(\alpha-\beta)\Delta}-1}\\&
=|F_{0}^{i}|^{2}+C'\Delta\left(\varepsilon a_{ii}p_i+\sum_{j\neq i}^{d}a_{ij}p_j+\sum_{j=1}^{d}b_{ij}p_je^{\beta\tau}\right)e^{\alpha\Delta}\frac{e^{(\alpha-\beta)k\Delta}-1}{e^{(\alpha-\beta)\Delta}-1}\\&
\quad-C'\Delta\varepsilon a_{ii}p_ie^{\alpha\Delta}\frac{e^{(\alpha-\beta)k\Delta}-1}{e^{(\alpha-\beta)\Delta}-1}.
	\end{aligned}$$

Since (\ref{changshu3}) holds for $\varepsilon<\frac{\min_i|a_{ii}|}{d\max_i|a_{ii}|}$, then there exist $\beta>0$ small enough and $\varepsilon'>0$ such that $\varepsilon\frac{\max_i|a_{ii}|}{\min_i|a_{ii}|}<\varepsilon'<\frac{1}{d}$ and
	\begin{equation}\label{changshu4}
		\varepsilon a_{ii}p_i+\sum_{j\neq i}a_{ij}p_j+\sum_{j=1}^{d}b_{ij}p_je^{\beta\tau}\leq-\varepsilon'\beta p_{i}
	\end{equation}
	for any $i \in \underline{d}$.

It follows that
$$\aligned e^{\alpha k\Delta}E|F_{k}^{i}|^{2}&\le |F_{0}^{i}|^{2}+ C'\Delta(-\varepsilon'\beta-\varepsilon a_{ii})p_ie^{\alpha\Delta}\frac{e^{(\alpha-\beta)k\Delta}-1}{e^{(\alpha-\beta)\Delta}-1}\\&
\le  |F_{0}^{i}|^{2}+C'\Delta\varepsilon'(\alpha-\beta)p_i\frac{e^{(\alpha-\beta)k\Delta}-1}{e^{(\alpha-\beta)\Delta}-1}.\endaligned$$
We have used in the last inequality the fact that $(-\varepsilon'\beta-\varepsilon a_{ii})e^{\alpha\Delta}\le \varepsilon'(\alpha-\beta)$ if $\alpha\in(\frac{\varepsilon}{\varepsilon'}\max_i|a_{ii}|,\min_i|a_{ii}|)$ and $\Delta$ is sufficiently small.

Moreover, by the fact that $x\le e^x-1, \forall x>0$, we have
$$\aligned e^{\alpha k\Delta}E|F_{k}^{i}|^{2}&< |F_{0}^{i}|^{2}+C'\varepsilon'p_i(e^{(\alpha-\beta)k\Delta}-1)\\&
=|F_{0}^{i}|^{2}-C'\varepsilon'p_i+C'\varepsilon'p_ie^{(\alpha-\beta)k\Delta}\\&
\le C'\varepsilon'p_ie^{(\alpha-\beta)k\Delta}.
\endaligned$$

Thus
\begin{equation}\label{fk}E|F_{k}^{i}|^{2}\le C'\varepsilon'p_ie^{-\beta k\Delta}.\end{equation}

On the other hand, since
$$\aligned|F_{k}^{i}|^{2}&=|X_{k}^{i}|^{2}-2\theta\Delta X_{k}^{i}f_k^i+\theta^2\Delta^2 |f_k^i|^2\\&
\ge |X_{k}^{i}|^{2}-\theta\Delta\left(\sum_{j=1}^{d}a_{ij}|X_{k}^{j}|^{2}+\sum_{j=1}^{d}b_{ij}|X_{k-[\frac{\tau(k\Delta)}{\Delta}]}^{j}|^{2}\right),\endaligned$$
then
$$(1-\theta\Delta a_{ii})|X_{k}^{i}|^{2}\le |F_{k}^{i}|^{2}+\theta\Delta\left(\sum_{j\neq i}a_{ij}|X_{k}^{j}|^{2}+\sum_{j=1}^{d}b_{ij}|X_{k-[\frac{\tau(k\Delta)}{\Delta}]}^{j}|^{2}\right).$$

Therefore
$$\aligned\sum_{i=1}^{d}\frac{|X_{k}^{i}|^{2}}{p_i}&\le\sum_{i=1}^{d}\frac{|F_{k}^{i}|^{2}}{p_i(1-\theta\Delta a_{ii})}+\sum_{i=1}^{d}\frac{\theta\Delta}{1-\theta\Delta a_{ii}}\sum_{j\neq i}\frac{a_{ij}p_j}{p_i}\frac{|X_{k}^{j}|^{2}}{p_j}\\&
\quad+\sum_{j=1}^d\left(\sum_{i=1}^{d}\frac{\theta\Delta}{1-\theta\Delta a_{ii}}\frac{b_{ij}p_j}{p_i}\right)\frac{|X_{k-[\frac{\tau(k\Delta)}{\Delta}]}^{j}|^{2}}{p_j}\\&
\le \sum_{i=1}^{d}\frac{|F_{k}^{i}|^{2}}{p_i(1-\theta\Delta a_{ii})}+\left(\sum_{i=1}^{d}\frac{\theta\Delta}{1-\theta\Delta a_{ii}}\max_{j\neq i}\frac{a_{ij}p_j}{p_i}\right)\sum_{j=1}^d\frac{|X_{k}^{j}|^{2}}{p_j}\\&
\quad+\sum_{j=1}^d\left(\sum_{i=1}^{d}\frac{\theta\Delta}{1-\theta\Delta a_{ii}}\frac{b_{ij}p_j}{p_i}\right)\frac{|X_{k-[\frac{\tau(k\Delta)}{\Delta}]}^{j}|^{2}}{p_j}.\endaligned$$

Then
$$\aligned\sum_{i=1}^{d}\frac{|X_{k}^{i}|^{2}}{p_i}&\le\frac{1}{1-D}\left(\sum_{i=1}^{d}\frac{|F_{k}^{i}|^{2}}{p_i(1-\theta\Delta a_{ii})}+\sum_{j=1}^d\left(\sum_{i=1}^{d}\frac{\theta\Delta}{1-\theta\Delta a_{ii}}\frac{b_{ij}p_j}{p_i}\right)\frac{|X_{k-[\frac{\tau(k\Delta)}{\Delta}]}^{j}|^{2}}{p_j}\right)\endaligned$$
where $D=\sum_{i=1}^{d}\frac{\theta\Delta}{1-\theta\Delta a_{ii}}\max_{j\neq i}\frac{a_{ij}p_j}{p_i}$.

By \eqref{fk} and inductive hypothesis \eqref{jiashe1}, it follows that	
$$\aligned\sum_{i=1}^{d}\frac{E|X_{k}^{i}|^{2}}{p_i}&\le\frac{1}{1-D}\left(\frac{dC'\varepsilon'e^{-\beta k\Delta}}{1-\theta\Delta \max_i|a_{ii}|}+C'e^{-\beta k\Delta}e^{\beta\tau}\sum_{j=1}^d\sum_{i=1}^{d}\frac{\theta\Delta}{1-\theta\Delta a_{ii}}\frac{b_{ij}p_j}{p_i}\right)\\&
=C'e^{-\beta k\Delta}\frac{1}{1-D}\left(\frac{d\varepsilon'}{1-\theta\Delta \max_i|a_{ii}|}+e^{\beta\tau}\sum_{j=1}^d\sum_{i=1}^{d}\frac{\theta\Delta}{1-\theta\Delta a_{ii}}\frac{b_{ij}p_j}{p_i}\right).\endaligned$$
Denote $G(\Delta):=\frac{1}{1-D}\left(\frac{d\varepsilon'}{1-\theta\Delta \max_i|a_{ii}|}+e^{\beta\tau}\sum_{j=1}^d\sum_{i=1}^{d}\frac{\theta\Delta}{1-\theta\Delta a_{ii}}\frac{b_{ij}p_j}{p_i}\right)$.

Since $G(0)=d\varepsilon'<1$, and $G$ is continuous, then there exists $\Delta^*$ small enough such that $G(\Delta)\le 1$ for all $\Delta\le \Delta^*.$

Then we have
$$\aligned\sum_{i=1}^{d}\frac{E|X_{k}^{i}|^{2}}{p_i}&\le C'e^{-\beta k\Delta}.\endaligned$$

Thus
$$E|X_{k}^{i}|^{2}\le C'p_ie^{-\beta k\Delta}, i=1,2,...,d.$$

That is, $\theta$-EM method is exponentially stable in mean square sense.

Almost surely exponential stability of $\theta$-EM method $X_k$ is a direct consequence of mean square exponential stability by using Chebyshev inequality and Borel-Cantelli Lemma (see e.g. \cite{1}).

We complete the proof. $\square$

\begin{Theorem}\label{num1}
Suppose all conditions in Theorem \ref{Theorem21} hold, and there exists a constant $K>0$ such that
\begin{equation}\label{xx}|f(x,y)|\le K(|x|+|y|).\end{equation}
Then for any fixed $0\le\theta\le\frac{1}{2},$ $\theta$-EM method is mean square and almost surely exponentially stable for sufficient small step size $\Delta>0$.
\end{Theorem}

\textbf{Proof.} Firstly, it still holds that
\begin{equation}
 \begin{aligned}
   |F^i_{k+1}|^{2}
   &\le |F^i_{k}|^{2}+[2X^i_{k}f_i(X_{k},X_{k-[\frac{\tau(k\Delta)}{\Delta}]})+ \sum_{l=1}^mg_{il}(X_{k},X_{k-[\frac{\tau(k\Delta)}{\Delta}]})^{2}\\
   &\quad+(1-2\theta)|f_i(X_{k},X_{k-[\frac{\tau(k\Delta)}{\Delta}]})|^{2}\Delta ]\Delta+M^{i}_{k}.
\end{aligned}
\end{equation}

Notice that linear growth condition (\ref{xx}) on $f$ implies
\begin{equation}
 \begin{aligned}|F^i_{k}|^{2}&=|X^i_{k}|^{2}-2\theta \Delta X^i_kf_i(X_{k},X_{k-[\frac{\tau(k\Delta)}{\Delta}]})+\theta^2\Delta^2|f_i(X_{k},X_{k-[\frac{\tau(k\Delta)}{\Delta}]})|^2\\&
 \le|X^i_{k}|^{2}+2\theta \Delta K|X^i_k|(|X_k|+|X_{k-[\frac{\tau(k\Delta)}{\Delta}]}|)+2\theta^2\Delta^2K^2(|X_k|^2+|X_{k-[\frac{\tau(k\Delta)}{\Delta}]}|^2)\\&
 \le(1+4\theta \Delta K+2\theta^2\Delta^2K^2)|X^i_{k}|^{2}+(\theta \Delta K+2\theta^2\Delta^2K^2)(\sum_{j\ne i}^d|X^j_k|^2+|X_{k-[\frac{\tau(k\Delta)}{\Delta}]}|^2),
 \end{aligned}
\end{equation}
then we have
$$|X^i_{k}|^{2}\geq \frac{|F^i_{k}| ^{2}-(\theta \Delta K+2\theta^2\Delta^2K^2)(\sum_{j\ne i}^d|X^j_k|^2+|X_{k-[\frac{\tau(k\Delta)}{\Delta}]}|^2)}{1+4\theta \Delta K+2\theta^2\Delta^2K^2}.$$

If $\theta\in[0,\frac{1}{2}],$ then by Assumption \ref{assumption 2.2} and (\ref{xx}), we have
\begin{equation}\label{dd}
 \begin{aligned}
   |F^i_{k+1}|^{2}
   &\le |F^i_{k}|^{2}+[\sum_{j=1}^{d}a_{i j}|X_{k}^j|^{2}+\sum_{j=1}^{d} b_{i j}|X_{k-[\frac{\tau(k\Delta)}{\Delta}]}^j|^{2}\\&\quad +2(1-2\theta)K^2(|X_{k}|^2+|X_{k-[\frac{\tau(k\Delta)}{\Delta}]}|^2)\Delta]\Delta+M^i_{k}\\&
   =|F^i_{k}|^{2}+(2(1-2\theta)K^2\Delta+a_{ii})|X^i_{k}|^2\Delta+\sum_{j\ne i}(2(1-2\theta)K^2\Delta+a_{ij})|X^j_{k}|^2\Delta\\&\quad+\sum_{j=1}^d(2(1-2\theta)K^2\Delta
   +b_{ij})|X^j_{k-[\frac{\tau(k\Delta)}{\Delta}]}|^2\Delta+M^i_{k}\\&\le |F^i_{k}|^{2}+(2(1-2\theta)K^2\Delta+a_{ii})\Delta\\&\quad\times\frac{|F^i_{k}| ^{2}-(\theta \Delta K+2\theta^2\Delta^2K^2)(\sum_{j\ne i}^d|X^j_k|^2+|X_{k-[\frac{\tau(k\Delta)}{\Delta}]}|^2)}{1+4\theta \Delta K+2\theta^2\Delta^2K^2}\\&\quad +\sum_{j\ne i}(2(1-2\theta)K^2\Delta+a_{ij})|X^j_{k}|^2\Delta\\&\quad+\sum_{j=1}^d(2(1-2\theta)K^2\Delta
   +b_{ij})|X^j_{k-[\frac{\tau(k\Delta)}{\Delta}]}|^2\Delta+M^i_{k}.
\end{aligned}
\end{equation}

We have used the fact that $2(1-2\theta)K^2\Delta+a_{ii}<0$ in the last inequality.

Thus
$$\aligned |F^i_{k+1}|^{2}&\le |F^i_{k}|^{2}\left(1+\frac{(2(1-2\theta)K^2\Delta+a_{ii})\Delta}{1+4\theta \Delta K+2\theta^2\Delta^2K^2}\right)\\&\quad+\sum_{j\ne i}\left(2(1-2\theta)K^2\Delta-(\theta K+2\theta^2\Delta K^2)\frac{(2(1-2\theta)K^2\Delta+a_{ii})\Delta}{1+4\theta \Delta K+2\theta^2\Delta^2K^2}+a_{ij}\right)|X^i_{k}|^2\Delta\\&\quad+\sum_{j=1}^d\left(2(1-2\theta)K^2\Delta-(\theta K+2\theta^2\Delta K^2)\frac{(2(1-2\theta)K^2\Delta+a_{ii})\Delta}{1+4\theta \Delta K+2\theta^2\Delta^2K^2}+b_{ij}\right)
\\&\quad\times|X^i_{k-[\frac{\tau(k\Delta)}{\Delta}]}|^2\Delta+M^i_{k}.\endaligned$$

It is obvious that $$C_{1,\Delta}:=\frac{2(1-2\theta)K^2\Delta+a_{ii}}{1+4\theta \Delta K+2\theta^2\Delta^2K^2}\to a_{ii},$$
$$C_{2,\Delta}:=2(1-2\theta)K^2\Delta-(\theta K+2\theta^2\Delta K^2)\frac{(2(1-2\theta)K^2\Delta+a_{ii})\Delta}{1+4\theta \Delta K+2\theta^2\Delta^2K^2}+a_{ij}\to a_{ij}$$
and
$$C_{3,\Delta}:=2(1-2\theta)K^2\Delta-(\theta K+2\theta^2\Delta K^2)\frac{(2(1-2\theta)K^2\Delta+a_{ii})\Delta}{1+4\theta \Delta K+2\theta^2\Delta^2K^2}+b_{ij}\to b_{ij}$$
as $\Delta\to0.$

Then for any $C<\min_i|a_{ii}|$, we can choose $\tilde{\varepsilon}>0$ and $\Delta>0$ small enough such that
$$|F^i_{k+1}|^{2}\le |F^i_{k}|^{2}(1-C\Delta)+\sum_{j\ne i}(a_{ij}+\tilde{\varepsilon})|X^j_{k}|^2\Delta+\sum_{j=1}^d(b_{ij}+\tilde{\varepsilon})|X^j_{k-[\frac{\tau(k\Delta)}{\Delta}]}|^2+M^i_{k}.$$

Note that (\ref{changshu3}) implies that there exist $\beta>0$ small enough and $\varepsilon'>0$ such that $\varepsilon\frac{\max_i|a_{ii}|}{\min_i|a_{ii}|}<\varepsilon'<\frac{1}{d}$ and
	\begin{equation}
		\varepsilon a_{ii}p_i+\sum_{j\neq i}(a_{ij}+\tilde{\varepsilon})p_j+\sum_{j=1}^{d}(b_{ij}+\tilde{\varepsilon})p_je^{\beta\tau}\leq-\varepsilon'\beta p_{i}
	\end{equation}
	for any $i \in \underline{d}$.

That is, (\ref{changshu4}) holds for $a_{ij}$ and $b_{ij}$ replaced by $a_{ij}+\tilde{\varepsilon}$ and $b_{ij}+\tilde{\varepsilon}$, respectively.

Then repeat the following part of the proof of Theorem \ref{Theorem21} from line to line, we complete the proof. $\square$

\section{Mean square exponential stability of MTEM method}

In this Section, we will prove that MTEM method replicates the mean exponential stability of exact solution under given conditions.

To obtain the mean exponential stability of MTEM method, we need the following two Lemmas.
	
	\begin{Lemma}\label{Lemma1}
		Suppose the Assumption \ref{assumption 2.1} and Assumption \ref{assumption 2.2} hold. Then for any fixed $\Delta>0$,there exist constant matrices $a_{i i}\in \mathbb{R};a_{i j}\ge 0,i\neq j ;b_{ij}\ge0.$ $i \in \underline{d}$, such that
		\begin{equation}\label{3.6}	
			2x_{i}\cdot f_{\Delta,i}( x, y)+\sum_{l=1}^{m} g^{2}_{\Delta,il}(x, y) \leq \sum_{j=1}^{d}a_{i j}x_{j}^{2}+\sum_{j=1}^{d} b_{i j}y_{j}^{2}, \quad i \in \underline{d}
		\end{equation}
		for any $x, y\in\mathbb{R}^d$.
	\end{Lemma}
	\textbf{Proof.}
	
	If $|x| \vee|y|  \leq h(\Delta)$, then \eqref{3.6} holds naturally by \eqref{2.4}.
	Now assume $|x| \vee|y|>h(\Delta)$. Then by $\eqref{2.4}$, we have
	$$
	\begin{aligned}
		&\quad2x_{i}f_{\Delta,i}( x, y)+\sum_{l=1}^{m}g_{\Delta,il}^{2}( x, y)\\
		&=2x_{i}\frac{|x| \vee|y|}{h(\Delta)} \cdot f_{i}(\frac{h(\Delta)}{|x| \vee|y|}(x,y))+\sum_{l=1}^{m}(\frac{|x| \vee|y|}{h(\Delta)} \cdot g_{il}(\frac{h(\Delta)}{|x| \vee|y|}(x,y)))^{2}\\
		& =(\frac{|x| \vee|y|}{h(\Delta)})^{2}\cdot[2x_{i}\frac{h(\Delta)}{|x| \vee|y|} \cdot f_{i}(\frac{h(\Delta)}{|x| \vee|y|}(x,y))+\sum_{l=1}^{m} g_{il}^{2}(\frac{h(\Delta)}{|x| \vee|y|}(x,y))]\\
		& \leq (\frac{|x| \vee|y|}{h(\Delta)})^{2}\cdot[\sum_{j=1}^{d}a_{i j}(x_{j}\frac{h(\Delta)}{|x| \vee|y|})^{2}+\sum_{j=1}^{d}b_{i j}(\frac{h(\Delta)}{|x| \vee|y|}x_{j})^{2}]\\
		&= \sum_{j=1}^{d}a_{i j}x_{j}^{2}+\sum_{j=1}^{d}b_{i j}y_{j}^{2}.
	\end{aligned}
	$$Then \eqref{3.6} holds.
	
	We complete the proof.

\begin{Lemma}\label{Lemma2}
		Suppose Assumption \ref{assumption 2.1} holds. Then for any fixed $\Delta>0,$ the modified truncated functions $f_\Delta$ and $g_\Delta$ are linear growing with coefficient $L_{h(\Delta)}$. That is
\begin{equation}\label{jieduan}	
			|f_\Delta(x,y)|\vee|g_\Delta(x,y)|\le L_{h(\Delta)}(|x|+|y|)
		\end{equation}
holds for all $x,y\in\mathbb{R}^d.$
	\end{Lemma}
\textbf{Proof.} Since $f_\Delta$ and $g_\Delta$ are defined in the same way, we only need prove linear growth condition of $f$. If $|x| \vee|y|  \leq h(\Delta)$, by definition of $f_\Delta$ and $g_\Delta$, it follows that
$$|f_\Delta(x,y)|=|f(x,y)|=|f(x,y)-f(0,0)|\le L_{h(\Delta)}(|x-0|+|y-0|).$$
If $|x|\vee|y|>h(\Delta)$, denote $a=\frac{h(\Delta)}{|x|\vee|y|}$. Since $|ax|\vee|ay|\le h(\Delta)$, then
$$|f_\Delta(x,y)|=\frac{1}{a}|f(ax,ay)|=\frac{1}{a}|f(ax,ay)-f(0,0)|\le \frac{1}{a}L_{h(\Delta)}(|ax-0|+|ay-0|).$$

We complete the proof.

\begin{Remark}
In \cite{1}, the author proved that both $f_\Delta$ and $g_\Delta$ are globally Lipschitz continuous for any fixed $\Delta>0$. However, we can only obtain $|f_\Delta(x,y)|\le 5L_{h(\Delta)}(|x|+|y|)$ while by Lemma \ref{Lemma2}, we have $|f_\Delta(x,y)|\le L_{h(\Delta)}(|x|+|y|)$.
\end{Remark}

Now we are ready to present the main result in this section.

	\begin{Theorem}\label{Theorem2}
		Suppose Assumptions \ref{assumption 2.1} and \ref{assumption 2.2} hold. Then the MTEM method \eqref{em} is mean quare and almost surely exponentially stable.
	\end{Theorem}
	\textbf{Proof.}
	By \eqref{em}, for any $k\geqslant 0$, we have
	\begin{equation}\label{3.7}
		\begin{aligned}
			&X_{k+1}^{i}=X_{k}^{i}+f_{\Delta,i}\left( X_{k}, X_{k-[\frac{\tau(k\Delta)}{\Delta}]}\right) \Delta
			+\sum_{l=1}^{m}g_{\Delta,il}\left( X_{k}, X_{k-[\frac{\tau(k\Delta)}{\Delta}]}\right) \Delta B^l_{k}.
		\end{aligned}
	\end{equation}
	So,
	\begin{equation}\label{3.8}
		\begin{aligned}
			|X_{k+1}^{i}|^{2}&=|X_{k}^{i}|^{2}+(2X_{k}^{i} f_{\Delta,i}\left( X_{k}, X_{k-[\frac{\tau(k\Delta)}{\Delta}]}\right)
			+\sum_{l=1}^{m}g_{\Delta,il}^{2}\left( X_{k}, X_{k-[\frac{\tau(k\Delta)}{\Delta}]}\right) )\Delta\\&\quad +f^{2}_{\Delta,i}\left( X_{k}, X_{k-[\frac{\tau(k\Delta)}{\Delta}]}\right) \Delta^{2} +M^i_{k}\\
		\end{aligned}
	\end{equation}
	where
	$$\begin{aligned}
		M^i_{k}&=2\langle X_{k}^{i},\sum_{l=1}^{m}g_{\Delta,il}\left( X_{k}, X_{k-[\frac{\tau(k\Delta)}{\Delta}]}\right)\Delta B^l_{k}\rangle\\&\quad+2\langle f_{\Delta,i}\left( X_{k}, X_{k-[\frac{\tau(k\Delta)}{\Delta}]}\right) ,\sum_{l=1}^{m}g_{\Delta,il}\left( X_{k}, X_{k-[\frac{\tau(k\Delta)}{\Delta}]}\right)\Delta B^l_{k}\rangle\Delta \\
		&\quad+\left(\left|\sum_{l=1}^{m}g_{\Delta,il}\left( X_{k}, X_{k-[\frac{\tau(k\Delta)}{\Delta}]}\right)\Delta B^l_{k}\right|^{2}-\sum_{l=1}^{m}g^2_{\Delta,il}\left( X_{k}, X_{k-[\frac{\tau(k\Delta)}{\Delta}]}\right)\Delta\right).
	\end{aligned}$$

	It is obvious that $M^i_{k}$ is a $\mathscr{F}_{k\Delta}$ martingale and $ EM^i_{k}=0$. Then, by  Lemmas \ref{Lemma1} and \ref{Lemma2}, we have
	$$	\begin{aligned}
		&\quad e^{\alpha(k+1)\Delta}E|X_{k+1}^{i}|^{2}-e^{\alpha k\Delta}E|X_{k}^{i}|^{2} \\
		&\leq e^{\alpha(k+1)\Delta}(1-e^{-\alpha  \Delta})E|X_{k}^{i}|^{2}+(\sum_{j=1}^{d}a_{ij}E|X_{k}^{j}|^{2}+\sum_{j=1}^{d}b_{ij}E|X_{k-[\frac{\tau(k\Delta)}{\Delta}]}^{j}|^{2})e^{\alpha(k+1)\Delta}\Delta\\
		&\quad+2L_{h(\Delta)}^{2}\Delta^{2}(E|X_{k}|^{2}+E|X_{k-[\frac{\tau(k\Delta)}{\Delta}]}|^{2})e^{\alpha(k+1)\Delta}\\
		&\leq e^{\alpha(k+1)\Delta}(1-e^{-\alpha  \Delta})E|X_{k}^{i}|^{2}+\sum_{j=1}^{d}(a_{ij}+2L_{h(\Delta)}^{2}\Delta)E|X_{k}^{j}|^{2}e^{\alpha(k+1)\Delta}\Delta\\
		&\quad+(\sum_{j=1}^{d} b_{ij}+2L^{2}_{h(\Delta)}\Delta)E|X^{j}_{k-[\frac{\tau(k\Delta)}{\Delta}]}|^{2}e^{\alpha(k+1)\Delta}\Delta.
	\end{aligned}$$

	Since (\ref{changshu}) holds, similar to the proof of Theorem \ref{Theorem1}, there exist $\beta>0$ and $\varepsilon>0$ small enough such that
	\begin{equation}\label{changshu1}
		\sum_{j=1}^{d}\left(a_{ij}+\varepsilon+e^{\beta\tau}(b_{ij}+\varepsilon)\right)p_{j}\leq-\beta p_{i}
	\end{equation}
	for any $i \in \underline{d}$.

Obviously, for the same $\bar{K}$, we have
	\begin{equation}\label{chushi}
		E|X_{k}^{i}|^{2}\leq \bar{K}p_{i}e^{-\beta k\Delta}E||\xi||^{2} ,\quad k=-\bar{m},-\bar{m}+1, \ldots, 0.
	\end{equation}

Now assume that
\begin{equation}\label{jiashe}
		E|X_{k}^{i}|^{2}\leq \bar{K}p_{i}e^{-\beta k\Delta}E||\xi||^{2} ,\quad  k\leqslant n-1, i=1,...,d
	\end{equation}
for some $\beta>0$.
	
Now for $k=n$, it follows that
$$\begin{aligned}
		e^{\alpha n \Delta}E|X_{n}^{i}|^{2}&\leq E|X_{0}^{i}|^{2}+\sum_{l=0}^{n-1} e^{\alpha(l+1)\Delta}(1-e^{-\alpha  \Delta})E|X_{l}^{i}|^{2}\\&\quad+\sum_{l=0}^{n-1}\sum_{j=1}^{d}(a_{ij}+2L_{h(\Delta)}^{2}\Delta)\Delta E|X_{l}^{j}|^{2}e^{\alpha(l+1)\Delta}\\
		&\quad+
		\sum_{l=0}^{n-1}\sum_{j=1}^{d}(b_{ij}+2L_{h(\Delta)}^{2}\Delta)\Delta E|X_{l-[\frac{\tau(l\Delta)}{\Delta}]}^{j}|^{2}e^{\alpha(l+1)\Delta}\\
		&= E|X_{0}^{i}|^{2}+\sum_{l=0}^{n-1} e^{\alpha(l+1)\Delta}(1-e^{-\alpha  \Delta}+(a_{ii}+2L_{h(\Delta)}^{2}\Delta)\Delta)E|X_{l}^{i}|^{2}\\
		&\quad+\sum_{l=0}^{n-1}\sum_{i\neq j}(a_{ij}+2L_{h(\Delta)}^{2}\Delta)\Delta e^{\alpha(l+1)\Delta}E|X_{l}^{j}|^{2}
		\\&\quad+\sum_{l=0}^{n-1}\sum_{j=1}^{d}(b_{ij}+2L_{h(\Delta)}^{2}\Delta)\Delta e^{\alpha(l+1)\Delta}E|X_{l-[\frac{\tau(l\Delta)}{\Delta}]}^{j}|^{2}.
\end{aligned}$$

If we choose $\alpha>\max\limits_{i\in \underline{d}}{|a_{ii}|}\vee\beta$, then it is obvious that $1-e^{-\alpha  \Delta}+(a_{ii}+2L_{h(\Delta)}^{2}\Delta)\Delta>0$ holds for any $\Delta>0.$ Therefore, by the induction hypothesis (\ref{jiashe}), we have
$$\begin{aligned}
e^{\alpha n \Delta}E|X_{n}^{i}|^{2}&\leq  E|X_{0}^{i}|^{2}+\sum_{l=0}^{n-1} e^{\alpha(l+1)\Delta}(1-e^{-\alpha  \Delta}+(a_{ii}+2L_{h(\Delta)}^{2}\Delta)\Delta)K_{1}p_{i}e^{-\beta l\Delta}\\
		&\quad+\sum_{l=0}^{n-1}\sum_{i\neq j}(a_{ij}+2L_{h(\Delta)}^{2}\Delta)\Delta e^{\alpha(l+1)\Delta}K_{1}p_{j}e^{-\beta l\Delta}
		\\&\quad+\sum_{l=0}^{n-1}\sum_{j=1}^{d}(b_{ij}+2L_{h(\Delta)}^{2}\Delta)\Delta e^{\alpha(l+1)\Delta}K_{1}p_{j}e^{-\beta (l-[\frac{\tau(l\Delta)}{\Delta}])\Delta}\\
		&\leq  E|X_{0}^{i}|^{2}+\sum_{l=0}^{n-1}(e^{\alpha\Delta}-1) K_{1}p_{i}\cdot e^{(\alpha-\beta)l\Delta}\\
		&\quad+\sum_{l=0}^{n-1}\sum_{j=1}^{d}((a_{ij}+2L_{h(\Delta)}^{2}\Delta)+(b_{ij}+2L_{h(\Delta)}^{2}\Delta)e^{\beta\tau})p_{j}\Delta  K_{1} e^{\alpha\Delta} e^{(\alpha-\beta)l\Delta}.
	\end{aligned}$$
where $K_1=\bar{K}E||\xi||^2.$
	
Then if we choose the stepsize $\Delta>0$ sufficiently small such that $2L_{h(\Delta)}^{2}\Delta\le\varepsilon$, \eqref{changshu1} yields
	$$	\begin{aligned}
		e^{\alpha n \Delta}E|X_{n}^{i}|^{2}
		&\leq  E|X_{0}^{i}|^{2}+K_{1}p_{i}\cdot(e^{\alpha\Delta}-1)\cdot \frac{e^{(\alpha-\beta)n\Delta}-1}{e^{(\alpha-\beta)\Delta}-1}-K_{1}\beta p_{i}\Delta e^{\alpha\Delta}\frac{e^{(\alpha-\beta)n\Delta}-1}{e^{(\alpha-\beta)\Delta}-1}\\
		&= E|X_{0}^{i}|^{2}+K_{1}p_{i}\cdot(e^{\alpha\Delta}-1-\beta\Delta e^{\alpha\Delta})\cdot \frac{e^{(\alpha-\beta)n\Delta}-1}{e^{(\alpha-\beta)\Delta}-1}.
	\end{aligned}$$

On the other hand, it is obvious that
	\begin{equation}\label{3.9}
		\begin{aligned}
			&e^{\alpha\Delta}-1-\beta\Delta e^{\alpha\Delta}=(1-\beta\Delta)e^{\alpha\Delta}-1
			\leq e^{-\beta\Delta}\cdot  e^{\alpha\Delta}-1=e^{(\alpha-\beta)\Delta}-1.
		\end{aligned}
	\end{equation}
	
Therefore,
	\begin{equation}\label{3.19}
		\begin{aligned}
			e^{\alpha n \Delta}E|X^{i}_{n}|^{2}&\leq E|X_{0}^{i}|^{2}+K_{1}p_{i}\cdot(e^{(\alpha-\beta)n\Delta}-1)\\
			&= E|X_{0}^{i}|^{2}-K_{1}p_{i}+K_{1}p_{i}\cdot e^{(\alpha-\beta)n\Delta} \\
			&\leq K_{1}p_{i}\cdot e^{(\alpha-\beta)n\Delta}.
		\end{aligned}
	\end{equation}
	
Thus,
	\begin{equation}\label{3.10}
		\begin{aligned}
			&E|X_{n}^{i}|^{2}\leq  K_{1}p_{i}\cdot e^{-\beta n \Delta}.
		\end{aligned}
	\end{equation}

Similar to almost sure exponential stability of $\theta$-EM method, a standard procedure of using Chebyshev inequality and Borel-Cantelli's Lemma implies that
$$\frac{\log |X_{n}^{i}|}{n\Delta}\le -\frac{\beta}{2}.$$
	We complete the proof.  $\square$

\section{Examples}
	
Now let us present two examples to interpret our theory.

	\textbf{Example 1}
	Consider a 2-D stochastic differential delay equation given by
	\begin{equation}\label{229}
		\begin{aligned}
			d x_{1}(t)&=  (-\frac{1}{5}x_{1}(t)-\frac{2}{9}x_{1}^{3}(t)+\frac{4}{5}x_{2}(t)+\frac{1}{10^{4}}x_{1}(t-\tau(t))+\frac{1}{10^{4}}x_{2}(t-\tau(t))) d t\\&\quad+\frac{2}{3}x^2_{1}(t) d B_{1}(t) \\
			dx_{2}(t)&=(\frac{\sqrt{38}}{625}x_{1}(t)-x_{2}(t)+\frac{1}{10^{4}} x_{1}(t-\tau(t))+\frac{1}{10^{4}}x_{2}(t-\tau(t)))d t+\sqrt{\frac{2}{5}}x_{2}(t)d B_{2}(t)\\
		\end{aligned}
	\end{equation}
	for $t \geq 0$, where the initial value $\xi(s)\in C([-\tau,0],\mathbb{R}^d), \tau(t)=0.1(1-|\sin(t)|)\le\tau=0.1$, and $B(\cdot)=\left(B_{1}(\cdot), B_{2}(\cdot)\right)^{T}$ is a 2-D Brownian motion.
	
Obviously, \eqref{229} is of the form \eqref{sdde} with $f(x, y):=\left(f_{1}(x, y), f_{2}(x, y)\right)^{T} \in \mathbb{R}^{2}$ defined by
	$$
	\begin{aligned}
		& f_{1}(x, y):=-\frac{1}{5}x_{1}-\frac{2}{9}x_{1}^{3}+\frac{4}{5}x_{2}+\frac{1}{10^{4}}y_{1}+\frac{1}{10^{4}}y_{2}, \\
		& f_{2}(x,y):=\frac{\sqrt{38}}{625}x_{1}-x_{2}+\frac{1}{10^{4}}y_{1}+\frac{1}{10^{4}}y_{2}\\
	\end{aligned}	
	$$
	for $x:=\left(x_{1}, x_{2}\right)^{T} \in \mathbb{R}^{2}, y:=\left(y_{1}, y_{2}\right)^{T} \in \mathbb{R}^{2}$, and
	$$
	g(x, y):=\left(\begin{array}{cc}
		\frac{2}{3}x^{2}_{1} & 0 \\
		0 & \sqrt{\frac{2}{5}}x_2
	\end{array}\right)
	$$
	for $x:=\left(x_{1}, x_{2}\right)^{T} \in \mathbb{R}^{2}, y:=\left(y_{1}, y_{2}\right)^{T} \in \mathbb{R}^{2}$. Clearly
	$$
	\begin{aligned}
		2x_{1} f_{1}(x, y)+g_{11}^{2}(x, y)+g_{12}^{2}(x, y)  &=-\frac{2}{5} x_{1}^{2}-\frac{4}{9}x_{1}x_{1}^{3}+\frac{8}{5}x_{1}x_{2}+\frac{2}{10^{4}}x_{1}y_{1}+\frac{2}{10^{4}}x_{1}y_{2}+\frac{4 }{9}x_{1}^{4} \\
		&\leq-\frac{399}{5000}x_{1}^{2}+2x_{2}^{2}+\frac{1}{10^{4}}y_{1}^{2}+\frac{1}{10^{4}}y_{2}^{2}\\
	\end{aligned}
	$$
	for any $x, y \in \mathbb{R}^{2}$ and
	$$
	\begin{aligned}
		2x_{2}f_{2}(x,y)+g_{21}^{2}(x,y)+g_{22}^{2}(x,y)  &= \frac{2\sqrt{38}}{625}x_{1}x_{2}-\frac{8}{5}x_{2}^{2}+\frac{2}{10^{4}}x_{2}y_{1}+\frac{2}{10^{4}}x_{2} y_{2}\\
		& \leq \frac{1}{5^{6}}x_{1}^{2}-\frac{399}{5000}x_{2}^{2}+\frac{1}{10^{4}}y_{1}^{2}+\frac{1}{10^{4}}y_{2}^{2}\\
	\end{aligned}
	$$
	for any $y\in\mathbb{R}^2$.

It is obvious that both $f$ and $g$ are locally Lipschitz continuous, and \eqref{3.1} holds for $K=2$. Thus, there exists a unique global solution to equation \eqref{229}. 	

On the other hand, if we take
	$$
	\begin{aligned}
		a_{11}=-\frac{399}{5000},a_{12}=2,a_{21}=\frac{1}{5^{6}},a_{22}=-\frac{399}{5000},
		b_{11}=\frac{1}{10^{4}},b_{12}=\frac{1}{10^{4}},b_{21}=\frac{1}{10^{4}},b_{22}=\frac{1}{10^{4}},
	\end{aligned}
	$$
	and choose $ \varepsilon=\frac{499}{1000}<\frac{1}{2}=\frac{\min_ia_{ii}}{d\max_ia_{ii}}, p_{1}=500, p_{2}=1$, we obtain that
	$$
	\begin{aligned}
		&\varepsilon a_{1 1}p_{1}+a_{1 2}p_{2}+b_{1 1}p_{1}+b_{1 2}p_{2}< 0\\
		&a_{2 1}p_{1}+\varepsilon a_{2 2}p_{2}+b_{2 1}p_{1}+b_{2 2}p_{2}< 0.\\
	\end{aligned}
	$$

That is, \eqref{w} holds with $a_{ij}, b_{ij}$ satisfies \eqref{changshu3} with given $p_i$ and $\varepsilon$, which clearly implies \eqref{changshu} for the same $a_{ij}, b_{ij}$ and $p_i$. Therefore, the trivial solution to \eqref{229} is exponentially stable in mean square sense by Theorem \ref{Theorem1}. Moreover, it is easy to verify that \eqref{danbian} holds for $L=1$. Thus $\theta$-EM  method is well defined in this case if $\Delta<1$. Then mean square and almost surely exponentially stable for the $\theta$-EM  method for any fixed $\theta\in(\frac{1}{2},1]$ by Theorem \ref{Theorem21}. Note that all assumptions in Theorem \ref{Theorem2} hold. Then Theorem \ref{Theorem2} implies that the corresponding MTEM method of \eqref{229} is also exponentially stable in mean square.
	
Meanwhile, we can find
	\begin{equation}\label{24.2}
		\begin{aligned}
			&2\langle x,f(x,y)\rangle+||g(x,y)||^2\\
			=&2x_{1} f_{1}(x, y)+2x_{2}f_{2}(x,y)+g_{11}^{2}(x, y)+g_{12}^{2}(x, y) +g_{21}^{2}(x,y)+g_{22}^{2}(x,y) \\
			=&-\frac{2}{5} x_{1}^{2}+(\frac{2\sqrt{38}}{625}+\frac{8}{5})x_{1}x_{2}-\frac{8}{5}x_{2}^{2}+\frac{2}{10^{4}}x_{1}y_{1}+\frac{2}{10^{4}}x_{1}y_{2}+\frac{2}{10^{4}}x_{2}y_{1}+\frac{2}{10^{4}}x_{2} y_{2}.\\
		\end{aligned}
	\end{equation}
We claim that \eqref{24.2} could not be written as Khasminskii-type condition \eqref{K}. Indeed, it follows that
	$$
	\begin{aligned}
		&2\langle x,f(x,y)\rangle+||g(x,y)||^2\\
		=&-\frac{2}{5} x_{1}^{2}+(\frac{2\sqrt{38}}{625}+\frac{8}{5})x_{1}x_{2}-\frac{8}{5}x_{2}^{2}+\frac{2}{10^{4}}x_{1}y_{1}+\frac{2}{10^{4}}x_{1}y_{2}+\frac{2}{10^{4}}x_{2}y_{1}+\frac{2}{10^{4}}x_{2} y_{2}\\
		\leq&-\frac{2}{5} x_{1}^{2}+(\frac{\sqrt{38}}{625}+\frac{4}{5})(n_{1}x_{1}^{2}+\frac{1}{n_{1}}x_{2}^{2})-\frac{8}{5}x_{2}^{2}+\frac{1}{10^{4}}(n_{2}x_{1}^{2}+\frac{1}{n_{2}}y_{1}^{2})+\frac{1}{10^{4}}(n_{3}x_{1}^{2}+\frac{1}{n_{3}}y_{2}^{2})\\
		&+\frac{1}{10^{4}}(n_{4}x_{2}^{2}+\frac{1}{n_{4}}y_{1}^{2})+\frac{1}{10^{4}}(n_{5}x_{2}^{2}+ \frac{1}{n_{5}}y_{2}^{2}) \\
		=&-\frac{1}{10^{4}}(4000-(16\sqrt{38}+8000)n_{1}-n_{2}-n_{3})x_{1}^{2}-\frac{1}{10^{4}}(16000-\frac{8000+16\sqrt{38}}{n_{1}}-n_{4}-n_{5})x_{2}^{2}\\
		&+\frac{1}{10^{4}}(\frac{1}{n_{2}}+ \frac{1}{n_{4}})y_{1}^{2}+\frac{1}{10^{4}}(\frac{1}{n_{3}}+ \frac{1}{n_{5}})y_{2}^{2}.
	\end{aligned}
	$$
	for any $n_{1},n_{2},n_{3},n_{4},n_{5}>0$.
	
However, if \eqref{K} holds for some $C_1>C_2>0$, then we must have
	
	\begin{equation}\label{24.3}
		\frac{1}{10^{4}}\min\left\lbrace 4000-(16\sqrt{38}+8000)n_{1}-n_{2}-n_{3},16000-\frac{8000+16\sqrt{38}}{n_{1}}-n_{4}-n_{5}\right\rbrace\ge C_1
	\end{equation}
	for some $n_{1},n_{2},n_{3},n_{4},n_{5}>0.$
	
Since there is no $n_1>0$ such that $4000-(16\sqrt{38}+8000)n_{1}>0$ and $16000-\frac{8000+16\sqrt{38}}{n_{1}}>0,$ \eqref{24.2} can never be written in the form of  Khasminskii-type conditions \eqref{K}. However, we can get mean square exponential stability of both the exact solution $x(t)$ and the $\theta$-EM method $X_k$ for any fixed $\theta\in(\frac{1}{2},1]$ by Theorem \ref{Theorem1} and Theorem \ref{Theorem21}. And $X_k$ is also almost surely exponentially stable.

\textbf{Example 2}
Consider the following 2-D stochastic differential delay equation
\begin{equation}\label{329}
	\begin{aligned}
		d x_{1}(t)&=(-\frac{2}{5}x_{1}(t)+\frac{4}{5}x_{2}(t)+\frac{1}{10^{4}}x_{1}(t-\tau(t))+\frac{1}{10^{4}}x_{2}(t-\tau(t))) d t\\&\quad+\frac{\sqrt{10}}{5}x_{1}(t) d B_{1}(t) \\
		dx_{2}(t)&=(\frac{\sqrt{38}}{625}x_{1}(t)-x_{2}(t)+\frac{1}{10^{4}} x_{1}(t-\tau(t))+\frac{1}{10^{4}}x_{2}(t-\tau(t)))d t\\&\quad+\frac{\sqrt{10}}{5}x_{2}(t) d B_{2}(t) \\
	\end{aligned}
\end{equation}
for $t \geq 0$, with initial value $\xi(s)\in C([-\tau,0],\mathbb{R}^d)$, where $\tau=0.1, \theta\in[-\tau,0],$ $\tau(t)=0.1|\sin(t)|\le\tau$.

Obviously, \eqref{329} is of the form \eqref{sdde} with $f(x, y):=\left(f_{1}(x, y), f_{2}(x, y)\right)^{T} \in \mathbb{R}^{2}$ defined by
$$
\begin{aligned}
	& f_{1}(x, y):=-\frac{2}{5}x_{1}+\frac{4}{5}x_{2}+\frac{1}{10^{4}}y_{1}+\frac{1}{10^{4}}y_{2}, \\
	& f_{2}(x,y):=\frac{\sqrt{38}}{625}x_{1}-x_{2}+\frac{1}{10^{4}}y_{1}+\frac{1}{10^{4}}y_{2}\\
\end{aligned}	
$$
for $x:=\left(x_{1}, x_{2}\right)^{T} \in \mathbb{R}^{2}, y:=\left(y_{1}, y_{2}\right)^{T} \in \mathbb{R}^{2}$, and
$$
g(x, y):=\left(\begin{array}{cc}
	\frac{\sqrt{10}}{5}x_{1} & 0 \\
	0 & \frac{\sqrt{10}}{5}x_{2}
\end{array}\right)
$$
for $x:=\left(x_{1}, x_{2}\right)^{T} \in \mathbb{R}^{2}, y:=\left(y_{1}, y_{2}\right)^{T} \in \mathbb{R}^{2}$. Similar to Example 1, it is easy to verify that \eqref{w} holds for
$$
\begin{aligned}
	a_{11}=-\frac{399}{5000},a_{12}=2,a_{21}=\frac{1}{5^{6}},a_{22}=-\frac{399}{5000},
	b_{11}=\frac{1}{10^{4}},b_{12}=\frac{1}{10^{4}},b_{21}=\frac{1}{10^{4}},b_{22}=\frac{1}{10^{4}},
\end{aligned}
$$
and \eqref{changshu3} holds for $ \varepsilon=\frac{499}{1000},p_{1}=500,p_{2}=1$. It is easy to see $f$ satisfies \eqref{xx}. Therefore,  it follows that the $\theta$-EM method (it is well defined since $f$ is global Lipschitz continuous in this case) is mean square and almost surely exponentially stable for any fixed $0\le\theta\le\frac{1}{2}$ by Theorem \ref{num1}.

Meanwhile, we can verify that \eqref{K} does not hold in the same way as in Example 1. However, we obtain that mean square exponential stability of both the exact solution $x(t)$ and the $\theta$-EM method $X_k$ for any fixed $0\le\theta\le\frac{1}{2}$ by Theorem \ref{Theorem1} and Theorem \ref{num1}. And $X_k$ is also almost surely exponentially stable.


\begin{thebibliography}{99}

\bibitem{GMY} Q. Guo, W. Liu, X. Mao, R. Yue, The truncated Milstein method for stochastic differential equations with commutative noise, J. Comput. Appl. Math., 338(2018)298-310.

    \bibitem{HW} E. Hairer, G. Wanner, Solving ordinary differential equation II: Stiff and differential-algebraic
problems, 2nd ed. Springer, Berlin, 1996.

\bibitem{1} G. Lan, Asymptotic exponential stability of modified truncated EM method for neutral stochastic differential delay equations, J. Comput. Appl. Math., 340(2018), 334-341.

\bibitem{LW} G. Lan, Q. Wang, Strong convergence rates of modified truncated EM methods for neutral stochastic differential delay equations, J. Comput. Appl. Math., 362(2019), 83-98.

\bibitem{LX} G. Lan, F. Xia, Strong convergence rates of modified truncated EM method for stochastic differential equations, J. Comput. Appl. Math., 334(2018), 1-17.

\bibitem{LMD} L. Liu, H. Mo, F. Deng, Split-step theta method for stochastic delay integro-differential equations with mean square exponential stability, Appl. Math. Comput., 353(2019) 320-328.

\bibitem{LZ} L. Liu, Q. Zhu, Mean square stability of two classes of theta method for neutral stochastic differential delay equations, J. Comput. Appl. Math., 305(2016) 55-67.

\bibitem{8} X. Mao, The truncated Euler-Maruyama method for stochastic differential equations, J. Comput. Appl. Math., 290(2015), 370-384.
		
		
\bibitem{9} X. Mao, Convergence rates of the truncated Euler-Maruyama method for stochastic differential equations, J. Comput. Appl. Math., 296(2016), 362-375.

\bibitem{14} X. Mao, Stochastic Differential Equations and Applications, second ed., Horwood, Chichester, 2007.

\bibitem{MY} X. Mao, C. Yuan, Stochastic differential equations with Markovian switching. Imperial College Press, London, 2006.

 \bibitem{NH} P. H. A. Ngoc, L. T. Hieu, A novel approach to mean square exponential stability of stochastic delay differential equations, IEEE Trans. Automat. Control, 66(2021), 2351-2356.

 \bibitem{OM} M. Obradovi$\acute{\textrm{c}}$, M. Milo$\check{\textrm{s}}$evi$\acute{\textrm{c}}$, Almost sure exponential stability of the $\theta $-Euler-Maruyama method, when $\theta \in (\frac{1}{2},1)$, for neutral stochastic differential equations with time-dependent delay under nonlinear growth conditions, Calcolo, 56:2(2019) 1-24.

 \bibitem{RAM} O.F. Rouz, D. Ahmadian, M. Milev, Exponential mean-square stability of two classes of theta Milstein methods for stochastic delay differential equations, AIP Conf. Proc., 1910, 060015 (2017); doi: 10.1063/1.5014009

\bibitem{WMS}  F. Wu, X. Mao, L. Szpruch, Almost sure exponential stability of numerical solutions for stochastic delay differential equations, Numer. Math. 115 (2010)681-697.

 \bibitem{ZSL} Y. Zhang, M. Song, M, Liu, Convergence and stability of stochastic theta method for nonlinear stochastic differential equations with piecewise continuous arguments, J. Comput. Appl. Math., 403(2022), 113849.

 \bibitem{ZYY} J. Zhao, Y. Yi, Y. Xu, Strong convergence and stability of the split-step theta method for highly nonlinear neutral stochastic delay integro differential equation, Appl. Numer. Math., 172(2022), 279-291.	
	\end{thebibliography}
\end{document}